\documentclass[graybox]{svmult}

\usepackage{mathptmx}       
\usepackage{helvet}         
\usepackage{courier}        
\usepackage{type1cm}        
\usepackage{makeidx}         
\usepackage{graphicx}        
\usepackage{multicol}        
\usepackage[bottom]{footmisc}

\usepackage{latexsym,amssymb,lastpage}
\usepackage{graphicx,amsfonts}
\usepackage{times,mathptmx,bm,amsmath}
\usepackage{dcolumn}

\usepackage{amsmath,amscd,amssymb,latexsym}
\usepackage{fullpage}
\usepackage{epsf}
\usepackage{epsfig}

\makeindex             

\newcommand{\beq}{\begin{equation}}
\newcommand{\eeq}{\end{equation}}

\def\vv{{\bf v}}

\def\u0{{\bf 0}}

\def\u0{{\bf 0}}
\def\b0{{\bf 0}}

\def\vE{ {{\bf E}} }
\def\vf{ {{\bf f}} }

\def\vT{ {{\bf T}} }
\def\vu{ {{\bf u}} }

\def\vv{ {{\bf v}} }

\def\vw{ {{\bf w}} }

\begin{document}

\title*{High Order Finite Difference Schemes for the Heat Equation Whose Convergence Rates are  Higher Than Their Truncation Errors$^{1,2}$}
\titlerunning{Schemes Whose Convergence Rates are  Higher Than Their Truncation Errors}
\author{A.~Ditkowski}
\institute{A. Ditkowski  \at School of Mathematical Sciences,
 Tel Aviv University,
 Tel Aviv 69978,
 Israel
\email{adid@post.tau.ac.il} \newline 
1. The final publication is available at link.springer.com via  http://dx.doi.org/10.1007/978-3-319-19800-2\_13\newline 
2.  In later works these schemes are called Error Inhibiting Schemes (EIS).
}
%
%
\maketitle


\vspace{-3cm}
\abstract{Typically when  a semi-discrete
approximation to a partial differential equation (PDE) is constructed  a discretization of the spatial operator with a truncation error
$\tau$ is derived. This discrete operator should  be semi-bounded for the scheme
to be stable. Under these conditions the Lax--Ricchtmyer equivalence
theorem assures that the scheme converges and that the error will
be, at most, of the order of $\| \tau \|$.  In most cases the error is in
indeed of the order of $\| \tau \|$. \newline \newline
 We demonstrate that for the Heat equation stable schemes can be constructed, whose   truncation errors are $\tau$, however, the actual errors are much smaller.
This gives more degrees of freedom in the design of schemes which
can make them more efficient (more accurate or compact) than
standard schemes. In some cases  the accuracy of the schemes can be
further enhanced using post-processing procedures.}

\section{Introduction}
\label{sec:1}

\setcounter{equation}{0} \setcounter{figure}{0}

Consider the differential problem:
\begin{eqnarray}\label{1.10}
\frac{\partial \, u}{\partial t} &= & P \left ( \frac{\partial \,
}{\partial x} \right ) u \;\;,\;\;\;\;\;  x\in\Omega\subset
{\mathbb{R}}^d \;,
t\ge0 \nonumber \\ 
 \\
u(t=0) &=& f(x) \;. \nonumber
\end{eqnarray}
where $u=(u_1,\ldots,u_m)^T$ and $P \left ( {\partial \,
}/{\partial x} \right ) $ is a linear  differential operator with appropriate boundary conditions. It is assumed that this problem is well posed, i.e. $\exists
K(t)<\infty$ such that  $|| u(t)|| \le K(t) || f ||$, where typically
$K(t)  = K e^{\alpha t}$.\newline

Let $Q$ be the discretization of $P \left ( \partial \, /\partial x
\right )$ where we assume:
\begin{description}
   \item[\textbf{Assumption 1}:] The discrete operator $Q$ is based on
  grid points $\{x_j \}$, $j=1, \ldots,N$.  \bigskip

  \item[\textbf{Assumption 2}:]  $Q$ is semi--bounded in some equivalent scalar
  product
$ \left( \cdot, \cdot \right )_H = \left( \cdot, H \cdot \right ) $,
i.e
\begin{equation}  \label{1.20}
\left( \vw, Q \vw \right )_H \,  \leq   \, \alpha \left( \vw, \vw
\right )_H \, = \, \alpha \left\| \vw   \right\|_H^2 \; .
\end{equation}

   \item[\textbf{Assumption 3}:] The local truncation error vector of $Q$
  is
  $\vT_e$ which is defined, at each entry $j$ by
   \begin{equation}  \label{1.30}
   \left(\vT_e \right)_j \, = \, \left( P w(x_j) \right)\, - \,  \left(Q \vw \right)_j ,
   \end{equation}
where $w(x)$ is a smooth function and $\vw$ is the projection of $w(x)$ onto the grid.  It is assumed that $ \left\| \vT_e \right\| \xrightarrow{N
\rightarrow \infty} 0$.
\end{description}

Consider the semi--discrete approximation:
\begin{eqnarray}\label{1.40} 
\frac{\partial \, \vv}{\partial t} &= & Q \vv \;\;,\;\;\;\;\;
t\ge0 
\\
\vv(t=0) &=& \vf \;. \nonumber
\end{eqnarray}
\bigskip
\textbf{Proposition}: Under Assumptions 1--3 The semi--discrete
approximation converges. \newline 
\textbf{Proof}: Let $\vu$ is
the projection of $u(x,t)$ onto the grid. Then, from assumption 3,
\begin{equation}\label{1.50} 
\frac{\partial \, \vu}{\partial t} \,=\, P \vu \,=\, Q \vu +  \vT_e
\; .
\end{equation}
Let $\vE = \vu - \vv$ then by subtracting \eqref{1.40} from
\eqref{1.50} one obtains the equation for $\vE$, namely
\begin{equation}\label{1.60} 
\frac{\partial \, \vE}{\partial t} \,=\,  Q \vE +  \vT_e\; .
\end{equation}
Next, by taking the $H$ scalar product with $\vE$, using
assumption 2 and the Schwartz inequality the following estimate can be derived
$$
\left(\vE,   \frac{\partial \, \vE}{\partial t} \right)_H \;= \;  \frac{1}{2}\frac{\partial \, }{\partial t}  \left(\vE,  \vE  \right)_H  \;=\;  \left\| \vE \right\|_H   \frac{\partial \, }{\partial t}\left| | \vE \right\|_H
\;=\;  \left(\vE,  Q \vE  \right)_H  + \left(\vE,  \vT_e  \right)_H  \; \nonumber \\
\; \leq\;   \alpha  \left(\vE, \vE  \right)_H \, + \, \left\| \vE
\right\|_H \left\| \vT_e \right\|_H \; . \nonumber
$$
Thus
\begin{equation}\label{1.80} 
\frac{\partial }{\partial t}  \left\| \vE \right\|_H \,\le \, \alpha
\left\| \vE \right\|_H +  \left\| \vT_e \right\|_H \; .
\end{equation}
Therefore:
\begin{equation}\label{1.90} 
\left\| \vE \right\|_H(t) \,\le \,\left\| \vE \right\|_H(t=0) {\rm
e}^{\alpha  t }  + \frac{{\rm e}^{\alpha  t }-1}{\alpha } \max_{0
\le\tau \le t} \left\| \vT_e \right\|_H \; \xrightarrow{N \rightarrow
\infty} \; 0\; .
\end{equation}
Here we assumed that $\left\| \vE \right\|_H(t=0)$ is either 0, or at
least of the order of machine accuracy. Equation \eqref{1.90}
establishes the fact that if the scheme is stable and consistent,
the numerical solution $\vv$ converges to the projection of the
exact solution onto the grid, $\vu$. Furthermore, it assures that
the error will be at most in the truncation error $\left\| \vT_e
\right\|_H$. This is one part of the landmark Lax--Richtmyer
equivalence theorem for semi--discrete approximation. See e.g.
\cite{gustafsson2013time}.

Due to this and similar results the common way for constructing
finite difference schemes is to derive a semi--bounded $Q$ with proper
truncation error. Typically the error $\left\| \vE \right\|_H$ is
indeed of the order of $\left\| \vT_e \right\|_H$.

It should be noted, however, that \eqref{1.60} is the exact equation
for the error dynamics, while \eqref{1.90} is an estimate. In this
paper we present finite difference schemes in which the errors are
smaller than their truncation errors. It is well known that boundary
conditions can be of one order lower accuracy without destroying the
convergence rate expected from the approximation at inner points,
see e.g. \cite{gustafsson1975convergence},
\cite{gustafsson1981convergence},   \cite{abarbanel2000error} and
\cite{svard2006order}. In \cite{gustafsson1981convergence} and
\cite{svard2006order} it was shown that for parabolic, incompletely
parabolic and 2nd-order hyperbolic equations the boundary conditions
can be  of two order less. Here, however, we consider low order
truncation errors in most or all of the grid points.

This paper is constructed as follows; in Section 2 we present a
preliminary example and illustrate the mechanism
which reduces the error. In Section 3 we present a two--point block
scheme which has a first order truncation error but has a second
or third order error. In Section 4 two three--point block schemes
are presented.
Discussion and remarks are presented in Section 5.


\section{Preliminary Example}
\label{sec:2}

Consider the heat equation
\begin{eqnarray}\label{2.10}
\frac{\partial \, u}{\partial t} &= &  \frac{\partial^2 \, u
 }{\partial x^2} +F(x,t) \;\;,\;\;\;\;\;  x\in [0, 2\pi)\;,
t\ge0 \nonumber
\\
u(t=0) &=& f(x) 
\end{eqnarray}
with periodic boundary conditions. \newline
 Let the scheme be:
\begin{eqnarray}\label{2.20}
\frac{\partial \, v_j }{\partial t}  &=& D_+D_- v_j \,+ \,(-1)^j
\, c \,v_j \,+ \, F_j(t)    \;\;;\;\; \;\;\; x_j=jh,\;\; h=2 \pi/N\;,\;\;N\;{\rm is\; even}\nonumber \\
  v_j(t=0) & = & f_j
\end{eqnarray}
where $ F_j(t) $ is the projection of $F(x,t)$ onto the spatial
 grid. The truncation error is
\begin{equation}\label{2.30} 
\left( T_e \right)_j = \frac{h^2}{12}\left( u_j \right )_{xxxx} +O(h^2) -
(-1)^j \, c \,v_j = O(1)  \; .
\end{equation}
 Formally, this scheme
is not consistent. The scheme \eqref{2.20} was run with the initial
condition $f(x) \, = \, \cos(x)$, $F=0$ and $N = 32, 64, \ldots,
1024$ with forward--Euler time propagator. The plots of $\log_{10}
\|\vE\|$ vs. $\log_{10} N$, at $t_{\rm final}=2 \pi$, for $c = 0,
0.5, 1$ are presented in Figure 1a. It can be seen that this is a
second order scheme.
\begin{figure}[h]
\begin{center}
\begin{tabular}{lll}
{\Large a}&&{\Large b}\\
\includegraphics[width=0.48\textwidth]{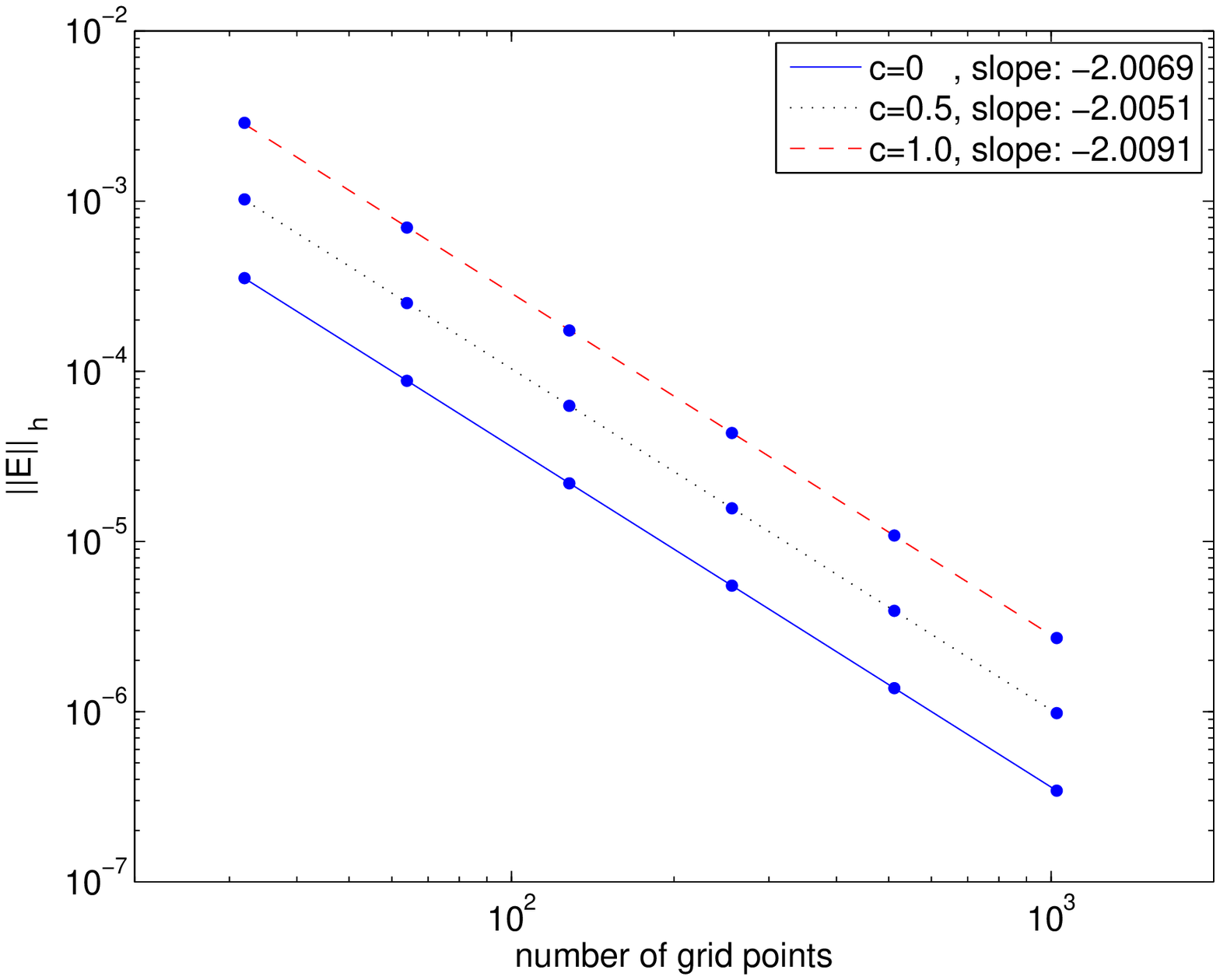}&&
\includegraphics[width=0.48\textwidth]{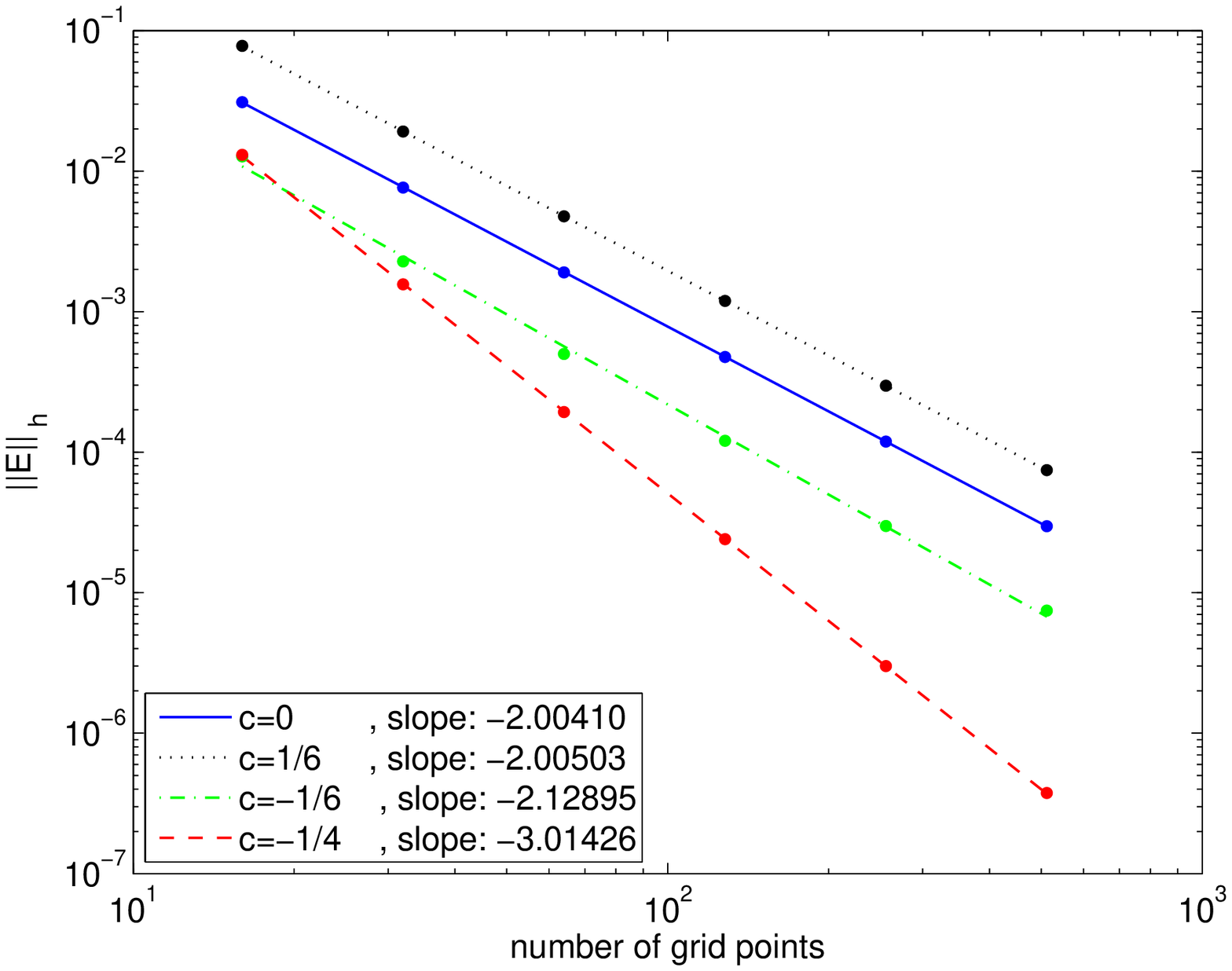}
\end{tabular}
\end{center}
\caption{Convergence plots, $\log_{10} \|\vE\|$ vs. $\log_{10} N$,
for different values of $c$. a: Scheme \eqref{2.20}. b: Scheme
\eqref{3.10}.
}
\label{fig:fig_1}
\end{figure}

In order to understand this phenomenon let us consider one high frequency mode of the error\footnote[3]{This scheme was presented for demonstrating the phenomenon that the error, due to high frequency modes, is lower than the truncation error. As this is not a practical scheme, full analysis of the error is not presented, only a   demonstration of the dynamics of high frequency error modes  is presented. Full analysis is given for the scheme presented in the next section}. Denote  $\left( T_c
\right)_j=(-1)^j \, c \,v_j$. This term can also be written as:
\begin{equation}\label{2.40} 
\left( T_c \right)_j \, = \, (-1)^j \, c \, h^{\alpha} \, = \, c
 h^{\alpha} \, {\rm e}^{i N x_j /2}
\end{equation}
For the scheme \eqref{2.20} $\alpha=0$. The equation for the error term caused by $T_c$ is
\begin{equation}\label{2.50} 
\frac{\partial \, \vE_c}{\partial t} \,=\,   D_+D_- \vE_c +  \vT_c
\end{equation}
Since $\left(  \vE_c \right)_j = \hat{\vE}_c \, {\rm e}^{i N x_j /2}\sqrt{2 \pi}$, $\hat{\vE}_c \in \mathbb{C}$,  then the equation for $\hat{\vE}_c$ is
\begin{equation}\label{2.60} 
\frac{\partial \, \hat{\vE}_c}{\partial t} \,=\,   -\left(
\frac{N}{2}\right)^2 \widehat{\vE}_c + c'
 h^{\alpha}
\end{equation}
Therefore,
\begin{eqnarray}\label{2.70}
\left\| {\vE}_c \right\|(t) \,=\,\left| \hat{\vE}_c \right|(t) &= & \left| \hat{\vE}_c \right|(0)
{\rm e}^{ -\left( \frac{N}{2}\right)^2 t } + \left(
\frac{2}{N}\right)^2 \left( 1 -{\rm e}^{ -\left(
\frac{N}{2}\right)^2 t } \right) c'  h^{\alpha}
 \\
& \le& \left| \hat{\vE}_c \right|(0) {\rm e}^{ -\left(
\frac{N}{2}\right)^2 t } + O\left( { h^{\alpha+2} } \right)\nonumber
\end{eqnarray}
Note that the actual error, $\left\| {\vE}_c \right\|(t)$, is two
orders lower that the truncation error, $\|T_c\|$. In the next
sections we present practical schemes which utilize this idea.

\section{Two--point Block, 3rd Order Scheme}

Let the grid be: $x_j = j \,h$, $h=2 \pi/(N+1)$ and $x_{j+1/2} =x_j
+h/2$, $j=0, \ldots N$. Altogether there are $2(N+1)$ points with
spacing of $h/2$. For simplicity, we assume that $N$ is even.

Consider the approximation:
\begin{eqnarray}\label{3.10}
\frac{d^2 }{dx^2} u_{j} &\approx & \frac{1}{(h/2)^2} \left[ \left(
u_{j-1/2} -2 u_{j} + u_{j+1/2}  \right) \, + \, c \left( -u_{j-1/2}
+3 u_{j} - 3 u_{j+1/2} +  u_{j+1} \right)  \right ]  \nonumber \\
\\
\frac{d^2 }{dx^2} u_{{j+1/2}} &\approx & \frac{1}{(h/2)^2} \left[
\left( u_{j} -2 u_{j+1/2} + u_{j+1}  \right) \, + \, c \left(
u_{j-1/2} -3 u_{j} + 3 u_{j+1/2} -  u_{j+1} \right)  \right ]\nonumber
\end{eqnarray}
The truncation errors are:
\begin{eqnarray}\label{3.20}
\left( T_e \right)_{j} \hspace{-0.2cm} & \; = \; &
\hspace{-0.2cm}\frac{1}{12} \left( \frac{h}{2}\right)^2 \left( u_j
\right )_{xxxx} +  c \left[ \left( \frac{h}{2}\right) \left( u_j
\right )_{xxx}  + \frac{1}{2} \left( \frac{h}{2}\right)^2 \left( u_j
\right
)_{xxxx}\right]+ O(h^3) = O(h)  \nonumber \\
\\
\left( T_e \right)_{j+\frac{1}{2}}\hspace{-0.2cm}  & \; = \; &
\hspace{-0.2cm}\frac{1}{12} \left( \frac{h}{2}\right)^2 \left(
u_{j+\frac{1}{2}}\right )_{xxxx} +  c \left[ -\left(
\frac{h}{2}\right) \left( u_{j+\frac{1}{2}} \right )_{xxx}  +
\frac{1}{2} \left( \frac{h}{2}\right)^2 \left(
u_{j+\frac{1}{2}}\right )_{xxxx}\right]+ O(h^3)  = O(h) \nonumber
\end{eqnarray}
The motivation leading to this scheme is that the highly oscillating
$O(h)$ error terms will be dissipated, as in the previous example,
while the $O(h^2)$ terms will be canceled, for the proper value of
$c$.

\subsection{ Analysis}

Let $\omega \in \{-N/2,\ldots,N/2$\} and
\begin{equation}\label{3.30} 
\nu \,=\,   \left\{ \begin{array}{lcll}
                       \omega -(N+1)& \hspace{2cm} & \omega >0\\
                       \\
                       \omega+(N+1)&                       & \omega \le0
                     \end{array}
\right.
\end{equation}
Then
\begin{equation}\label{3.40} 
 {\rm e}^{i \omega x_j} = {\rm e}^{i \nu x_j}\;\; {\rm and } \;\; {\rm e}^{i \omega
x_{j+1/2}} = -{\rm e}^{i \nu x_{j+1/2}} \;.
\end{equation}
We look for eigenvectors in the form of:
\begin{equation}\label{3.50} 
\psi_k(\omega) = \frac{\alpha_k}{\sqrt{2 \pi}} \left(
                          \begin{array}{c}
                            \vdots \\
                            {\rm e}^{i \omega x_j} \\
                            {\rm e}^{i \omega x_{j+1/2}} \\
                            \vdots
                          \end{array}
                        \right) +\frac{\beta_k}{\sqrt{2 \pi}} \left(
                          \begin{array}{c}
                            \vdots \\
                            {\rm e}^{i \nu x_j} \\
                            {\rm e}^{i \nu x_{j+1/2}} \\
                            \vdots
                          \end{array}
                        \right)
\end{equation}
where, for normalization, it is require that
$|\alpha_k|^2+|\beta_k|^2=1 $, $k=1,2$. The  expressions for
$\alpha_k$, $\beta_k$ and the eigenvalues (symbols) $\hat{Q}_k$ are:

\begin{equation}\label{3.52} 
\alpha_1 = \left[ {\sqrt{1+\frac{c^2 \cos (4 (h/2) \omega )+4 (2
c-1) \Delta
   \cos ((h/2) \omega )+4 (c (7 c-8)+2) \cos (2 (h/2) \omega )+(35
   c-32) c+8}{2 c^2 (2 \sin ((h/2) \omega )+\sin (2 (h/2) \omega
   ))^2}}} \; \right ]^{-1}
\end{equation}

\begin{equation}\label{3.53} 
\beta_1 = -\frac{i ((8 c-4) \cos ((h/2) \omega )+\Delta )}{2 c (2
\sin
   ((h/2) \omega )+\sin (2 (h/2) \omega )) \alpha_1^{-1} }
\end{equation}

\begin{equation}\label{3.54} 
\beta_2 = \left [ {\sqrt{1+\frac{2 c^2 (2 \sin ((h/2) \omega )+\sin
(2 (h/2)
   \omega ))^2}{c^2 \cos (4 (h/2) \omega )+4 (1-2 c) \Delta  \cos
   ((h/2) \omega )+4 (c (7 c-8)+2) \cos (2 (h/2) \omega )+(35
   c-32) c+8}}} \; \right ]^{-1}
\end{equation}

\begin{equation}\label{3.55} 
\alpha_2 = -\frac{2 i c (2 \sin ((h/2) \omega )+\sin (2 (h/2) \omega
   ))}{((4-8 c) \cos ((h/2) \omega )+\Delta ) \beta_2^{-1} }
\end{equation}
where

\begin{equation}\label{3.56} 
\Delta   = \sqrt{2 c^2 \cos (4 (h/2) \omega )+38 c^2+8 (c-1) (3 c-1)
\cos (2 (h/2) \omega )-32 c+8}
\end{equation}
and the

\begin{equation}\label{3.58} 
\hat{Q}_{1,2}(\omega)  = \frac{-4 +2 c (\cos (2 (h/2) \omega
   )+3) \pm  \Delta   }{2 (h/2)^2} \;.
\end{equation}
It can be shown that the eigenvalues are real and non positive for
$|c|<1/2$. It should be noted that by construction $\psi_k(\omega) \perp \psi_k(\nu)$, $\omega \ne \nu$. However, for all $ \omega$ $\psi_1(\omega)$ and $\psi_2(\omega)$ are not orthogonal. Nevertheless it can be shown that they are 'almost' perpendicular in the sense that 
$$
\cos\left( \theta(\omega)\right )\; = \; \frac{ \left \langle \psi_1(\omega),\psi_2(\omega) \right\rangle}{\|\psi_1(\omega)\| \; \|\psi_2(\omega)\|} \; > \; 0.9
$$
Therefore the scheme is stable.

For $\omega h \ll 1$ the eigenvalues and eigenvectors are:
\begin{equation}\label{3.60} 
\hat{Q}_1(\omega) = -\omega ^2 +\frac{(1+4 c)  \omega
   ^4}{12-24 c} \left( \frac{h}{2} \right)^2 + O(h^4)
\end{equation}
\begin{equation}\label{3.70} 
\alpha_1 = 1-\frac{c^2 }{32 (1-2
   c)^2} \left(\frac{\omega h}{2} \right)^6+O\left(\text{h}^7\right)\;, \;\;\; \beta_1 = -\frac{i c  }{4 - 8 c}\left(
\frac{\omega h}{2} \right)^3 + O(h^5)
\end{equation}
and
\begin{equation}\label{3.80} 
\hat{Q}_2(\omega) = -\frac{4 -8 c}{\text{(h/2)}^2}+(1-4 c) \omega ^2
+ O(h^2)
\end{equation}
\begin{equation}\label{3.90} 
\alpha_2 =  \frac{i c  }{2 c-1}\left( \frac{\omega h}{2} \right) +
O(h^3)\;,\;\;\;\beta_2 = 1 + O(h^2)
\end{equation}

If the initial condition is
\begin{equation}\label{3.100} 
\vv_j(0) = {\rm e}^{i \omega x_j}\;, \vv_{j+\frac{1}{2}}(0) = {\rm
e}^{i \omega x_{j+\frac{1}{2}}} \;\; ; \; \;\;\; \omega^2 h\ll 1
\end{equation}
then
\begin{eqnarray}\label{3.110}
\left(\vv\right)_j(t) &=& {\rm e}^{-\omega^2 t} \left(1 {-
\frac{(1+4 c) \omega^2 t  }{12-24 c}\left( \frac{\omega h}{2}
\right)^2 + O(h^4)}\right ){\rm e}^{i \omega x_j}+  { \left (
-\frac{i c }{4-8 c}\left( \frac{\omega h}{2} \right)^3 + O(h^5)
\right ){\rm e}^{i \nu x_j} }
\end{eqnarray}

The same expression hold for $x_{j+\frac{1}{2}}$. Therefore the
scheme is, in general, 2nd order and it is 3rd order if $c=-1/4$.
Note that by naive analysis of the truncation error terms,
\eqref{3.20}, one would expect to get 3rd order with $c=-1/6$.

We used the approximation \eqref{3.10} for solving the heat equation
\eqref{2.10}, where $F(x,t)$ and the initial condition were chosen
such that the exact solution is $u(x,t)=\exp(\cos(x-t))$. The scheme
was run with $N = 32, 64, \ldots, 1024$. 4th order Runge--Kutta
 scheme was used for time integration. The plots of $\log_{10} \|\vE\|$ vs.
$\log_{10} N$ for $c = 0, 1/6, -1/6,-1/4$ are presented in Figure
1b. As can be seen, the results are as predicted by the analysis.


\section{Three--point Block}

In this section we briefly present two schemes which are based on
three--point block.  Here we use  the grid, $x_j = j \,h$, $h=2
\pi/(N+1)$ with the internal block nodes $x_{j+1/3} =x_j +h/3$ and
$x_{j+2/3} =x_j +2 h/3$, $j=0, \ldots N$. Altogether there are
$3(N+1)$ points with spacing of $h/3$.

\bigskip
{\bf three--point block, 3rd order scheme}

 Consider the approximation:
\begin{eqnarray}\label{4.10}
\frac{d^2 }{dx^2} u_{j} &= & \frac{1}{4 (h/3)^2} \left[ \left( 4
u_{j-1/3} -8 u_{j} + 4 u_{j+1/3}  \right) \, + \, c \left(
-u_{j-1/3}
+3 u_{j} - 3 u_{j+1/3} +  u_{j+2/3} \right)  \right ] \, + \, O(h)\nonumber \\
\nonumber \\
\frac{d^2 }{dx^2} u_{j+1/3} &= & \frac{1}{4 (h/3)^2} \left[ \left(
4 u_{j} -8 u_{j+1/3} + 4 u_{j+2/3}  \right)  \right ] \, + \, O(h^2)\\
\nonumber \\
\frac{d^2 }{dx^2} u_{{j+2/3}} &= & \frac{1}{4 (h/3)^2} \left[ \left(
4 u_{j+1/3} -8 u_{j+2/3} + 4 u_{j+1}  \right) \, + \, c \left( u_{j}
-3 u_{j+1/3} + 3 u_{j+2/3} -  u_{j+1} \right)  \right ]\, + \, O(h)
\nonumber
\end{eqnarray}

This scheme was run under the same conditions as the example  in
Section 3. As in the previous example the truncation error is of
$O(h)$, however this is a 2nd order scheme and 3rd order for
$c=1.340$. The convergence results are presented in Figure 2a.

\begin{figure}[h]
\begin{center}
\begin{tabular}{lllll}
{\Large a}&&{\Large b}\\
\includegraphics[width=0.48\textwidth]{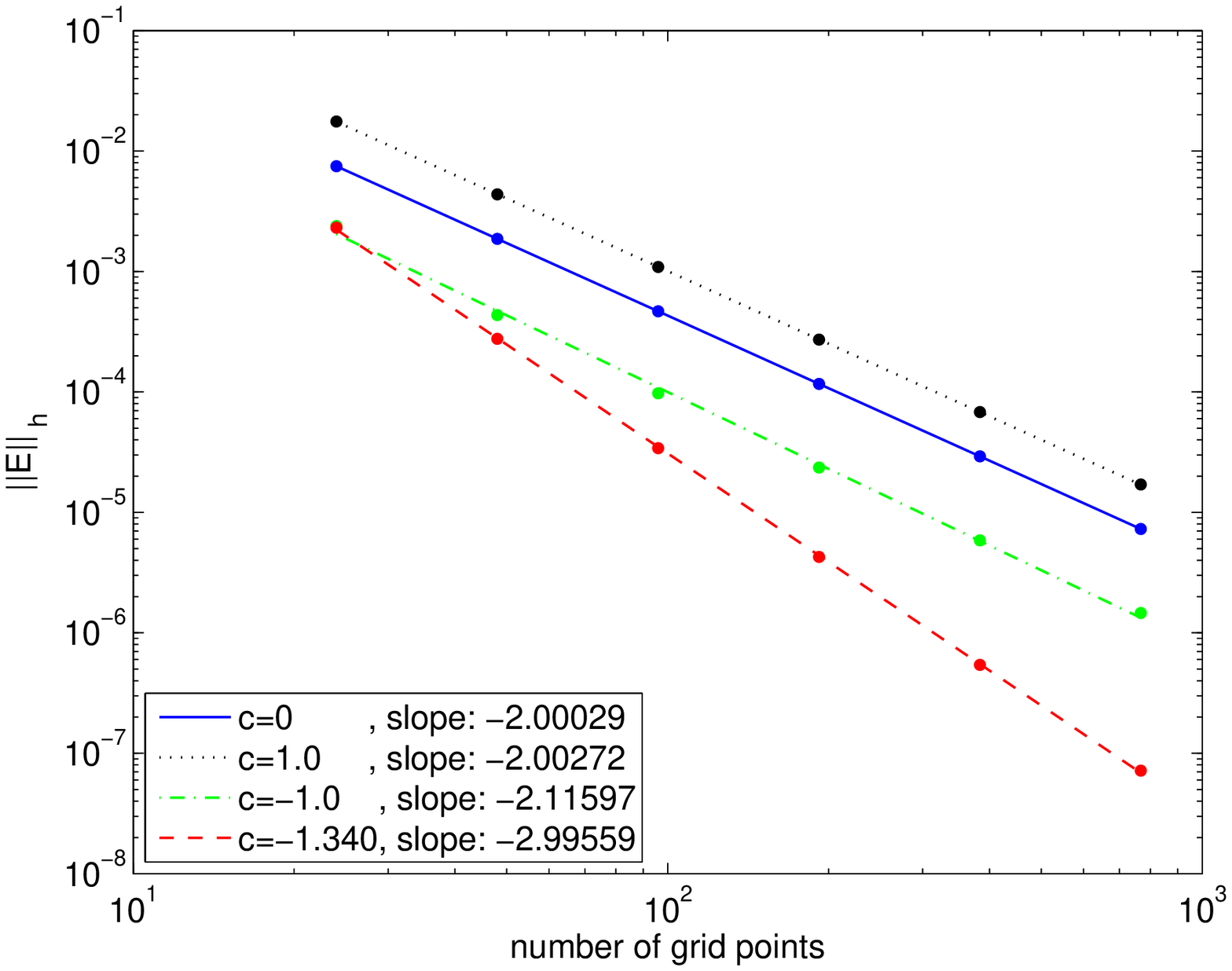}&\hspace{2cm}&
\includegraphics[width=0.48\textwidth]{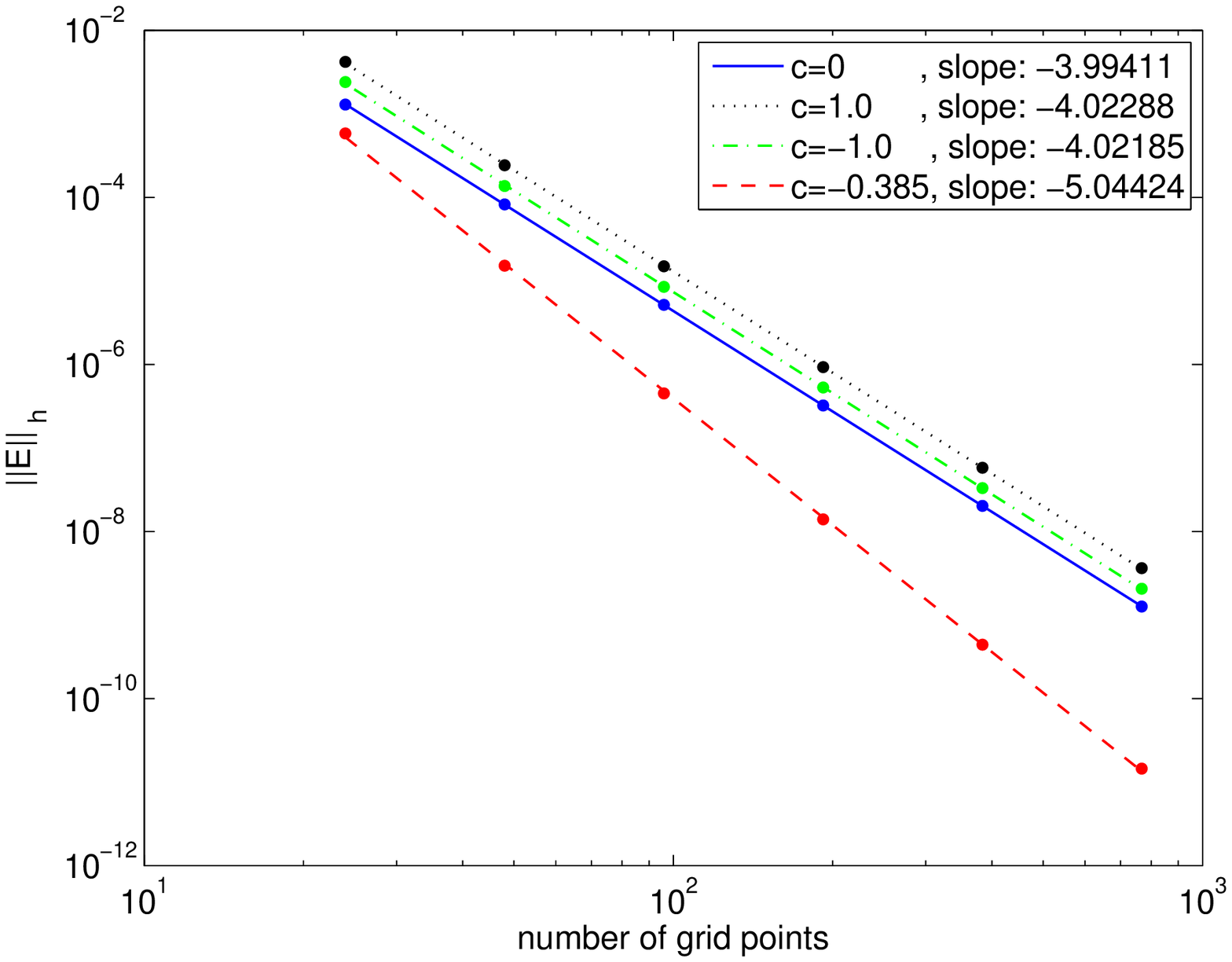}
\end{tabular}
\end{center}
\caption{Convergence plots, $\log_{10} \|\vE\|$ vs. $\log_{10} N$,
for different values of $c$. a: Scheme \eqref{4.10}. b: Scheme
\eqref{4.20}.
} \label{fig:fig_2}
\end{figure}

%
%
\bigskip
\newpage
{\bf three--point block, 5th order scheme}

\begin{eqnarray}\label{4.20}
\frac{d^2 }{dx^2} u_{j} &= & \frac{1}{12 (h/3)^2} \left[ \left(
-u_{j-2/3} + 16 u_{j-1/3} -30 u_{j} + 16 u_{j+1/3} -u_{j+2/3}
\right) \, + \, \right . \nonumber \\
&& \left . c \left( u_{j-2/3}-5u_{j-1/3}
+10 u_{j} - 10 u_{j+1/3} +  5u_{j+2/3}-u_{j} \right)  \right ] \, + \, O(h^3)\nonumber \\
\nonumber \\
\frac{d^2 }{dx^2} u_{j+1/3} &= & \frac{1}{12 (h/3)^2} \left[ \left(
-u_{j-1/3}+16 u_{j} -30 u_{j+1/3} + 16 u_{j+2/3} -u_{j}  \right)  \right ] \, + \, O(h^4)\\
\nonumber \\
\frac{d^2 }{dx^2} u_{{j+2/3}} &= & \frac{1}{12 (h/3)^2} \left[
\left( -u_{j}+16 u_{j+1/3} -30 u_{j+2/3} + 16 u_{j+1} -u_{j+4/3}
\right) \, + \, \right . \nonumber \\
&& \left . c \left( -u_{j-1/3} +5 u_{j} -10 u_{j+1/3} + 10 u_{j+2/3}
-  5 u_{j+1} + u_{j+4/3}\right) \right ] \, + \, O(h^3) \nonumber
\end{eqnarray}
This scheme was run under the same conditions as the previous
examples, with the exception that now  a 6th order Runge--Kutta scheme
was used for time integration. Here the truncation error is of
$O(h^3)$, however this is a 4th order scheme and 5th order for
$c=-0.385$. The convergence results are presented in Figure 2b.

 It should be noted that by taking $c=1$ the
coefficients of $u_{j-2/3}$ and $u_{j+4/3}$ are 0. Therefore the
scheme is more 'compact' than standard explicit 4th order scheme in
the sense that the scheme depends only on one term, on each side,
out of the three--point block $u_j, u_{j+1/3}, u_{j+2/3}$.
Potentially, this thinner stencil helps in the derivation of
boundary schemes for initial--boundary value problems.


\section{Summary}

In this paper we presented a few block--finite--difference schemes
in which the actual errors are much smaller than their truncation
errors. This reduction of error was achieved  by constructing the
truncation errors to be oscillatory and using the dissipative
property of the scheme.

A comparison between standard and block  finite difference schemes
in terms of the number of points out of the cell and operation count
is presented in Table 1. As can be seen, in terms of accuracy and
computational cost the 3rd and 5th order schemes are between the
standards 2nd and 4th order  and 4th and 6th order schemes
respectively. 
\begin{table}[h]
\begin{tabular}{ccccc|ccccccccccccc}
Scheme & Points out the block & \hspace{2cm}&\multicolumn{2}{c}{Number of operations} & Scheme & Points out the block & \hspace{2cm}&\multicolumn{2}{c}{Number of operations} \\
            &  (at each side) &  &  $+$& $\times$     &    &  (at each side) &  &  $+$& $\times$       \\
\hline
standard 2nd order& 1& & 2 & 3 &2--point block   3rd order& 1& &  3& 4 &   \\
standard 4th order& 2& & 4 & 5& 3--point block   3rd order& 1& & $2 \frac{2}{3}$ &  $3 \frac{2}{3}$\\
standard 6th order& 3&&  6 & 7 & 3--point block   5th order& 2& & $4 \frac{2}{3}$ &  $5 \frac{2}{3}$\\
&&&&& 3--point block   4th order& 1& & 4 & 5 \\
&&&&& compact
%
\end{tabular}
\caption{Comparison between standard finite difference schemes and
the ones presented on Sections 3,4.}
\end{table}

If  $c=-1/4$ in the two--point block, 3rd order
scheme \eqref{3.10}, the leading term in the error is highly
oscillatory, see \eqref{3.110}. It was suggested by Jennifer K. Ryan
\cite{J_R_pc} that this term can be filtered by post--processing. In
this technique the high frequency   error terms are filtered using
convolution with a proper kernel. This method was successfully
applied to  the discontinuous Galerkin method. As this filtering is
done only once, after the final time step, the cost is minimal, see
e.g. \cite{ryan2005extension}, \cite{cockburn2000post},
\cite{cockburn2003enhanced}. Here we used a global spectral filter
and the local filter suggested in \cite{cockburn2003enhanced}. As
can be seen in figure 3, the filtered scheme  is 4th order accurate.
The difference between both kernels is that the global filter is
more computationally expensive, $O(N \, \log N)$ for proper values
of $N$, but is accurate for small values of $N$ while the local
filter requires only $O(N)$ operation but is accurate only for large
values of $N$. The investigation of which is the optimal filter for
these schemes is a topic for future research.
\begin{figure}[h]
\begin{center}
\begin{tabular}{lll}
\includegraphics[width=0.48\textwidth]{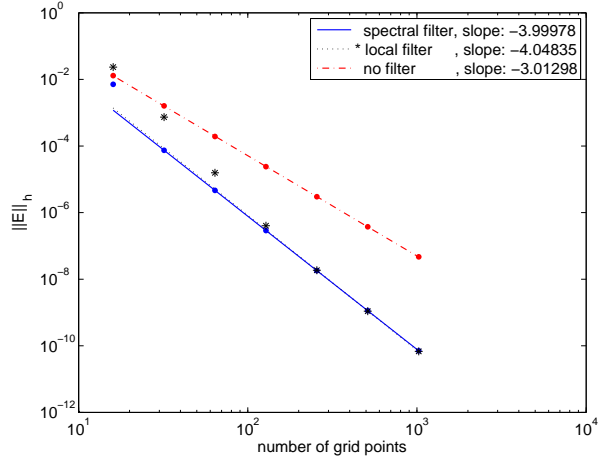}&&
\end{tabular}
\end{center}
\caption{Convergence plots, $\log_{10} \|\vE\|$ vs. $\log_{10} N$,
for scheme \eqref{3.10},  $c-1/4$. with no post--processing, with
spectral filter and with the filter suggested in
\cite{cockburn2003enhanced}.
} \label{fig:fig_3}
\end{figure}

 It should be noted that that finite difference
representations of discontinuous Galerkin schemes for the heat
equation have a similar form to the schemes presented above. They
may also present similar enhancing of accuracy. See
\cite{zhang2003analysis}. This manuscript was the inspiration to the
current work. As was also pointed out in \cite{zhang2003analysis}
the increasing of accuracy may be related to a phenomenon called
Supra--Convergence \cite{kreiss1986supra}.

Further research of the properties and implementations of these
schemes as well as the existence of these kinds of schemes for
hyperbolic problems are topics for future research.

\bigskip

{ \bf Acknowledgments} \newline The author would like to thank
Jennifer K. Ryan, Chi--Wang Shu and Sigal Gottlieb
 for the fruitful discussions and their help.
The author would also like to thank the anonymous reviewers for
their useful  remarks.



\vspace{-0.5cm}
\bibliographystyle{amsplain}

\bibliography{myref}

\end{document}